\documentclass[12pt,reqno,a4paper,notitlepage,oneside]{article}
\usepackage{amssymb,amsmath,amsthm,amscd} %пакеты AMS

\newcommand{\SJ}{\mathrm{SJ}}
\newcommand{\Jord}{\mathrm{Jord}}
\newcommand{\As}{\mathrm{As}}
\newcommand{\Alg}{\mathrm{Alg}}
\newcommand{\SId}{\mathrm{SId}}
\newcommand{\Di}{\mathrm{Di}}
\newcommand{\Cur}{\mathop{\mathrm{Cur}}\nolimits}
\newcommand{\Hom}{\mathcal{H}}
\newcommand{\Herm}{\mathrm{H}}
\newcommand{\Var}{\mathrm{Var}}
\newcommand{\charact}{\mathop{\mathrm{char}}\nolimits}
\newcommand{\Span}{\mathop{\mathrm{Span}}\nolimits}
\newcommand{\Coeff}{\mathop{\mathrm{Coeff}}\nolimits}

\newcommand{\mbdv}{\mathbin\dashv}
\newcommand{\mbvd}{\mathbin\vdash}
\newcommand{\mbdvjord}{\mathbin{{}_{(\dashv)}}}
\newcommand{\mbvdjord}{\mathbin{{}_{(\vdash)}}}
\newcommand{\dvjord}{{}_{(\dashv)}}
\newcommand{\vdjord}{{}_{(\vdash)}}
\newcommand{\co}[1]{\mathbin{{}_{(#1)}}}

\theoremstyle{plain}
\newtheorem{thm}{Theorem}
\newtheorem{lemma}[thm]{Lemma}
\newtheorem{prop}[thm]{Proposition}
\newtheorem{cor}[thm]{Corollary}

\theoremstyle{remark}
\newtheorem*{rem}{Remark}

\usepackage{hyperref}

\begin{document}

\linespread{1.25}\selectfont

\title{Special and exceptional Jordan dialgebras}

\author{Vasily~Voronin\thanks{The author was partially supported by
ADTP (Grant 2.1.1.10726), RFBR (Grant 09-01-00157-A), SSc-3669.2010.1, and Federal
Aim Program (contracts N. 02.740.11.5191, N. 14.740.11.0346).}}
%\address{Novosibirsk State University, Novosibirsk, Russia}
%\email{voronin.vasily@gmail.com}

\date{}

\maketitle

\begin{abstract}
In this paper, we study the class of Jordan dialgebras (also called quasi-Jordan
algebras). We develop an
approach for reducing problems on dialgebras to the case of ordinary
algebras. It is shown that straightforward generalizations of the
classical Cohn's, Shirshov's, and Macdonald's Theorems do not hold for
dialgebras. However, we prove dialgebraic analogues of these
statements. Also, we study multilinear special identities which
hold in all special Jordan algebras and do not hold in all Jordan algebras.
We find a natural correspondence between special identities for
ordinary algebras and dialgebras.
\end{abstract}

\section*{INTRODUCTION}

One of the most important classes of nonassociative algebras is the
class of Lie algebras defined by the anti-commutativity and Jacobi
identities $x^2=0$, $(xy)z+(zx)y+(yz)x=0$. It is well-known that
every associative algebra $A$ turns into a Lie algebra with respect
to the new product
 $[a,b]=ab-ba$, $a,b\in A$.
The Lie algebra obtained is denoted by $A^{(-)}$. The classical
Poincar\'e---Birkhoff---Witt Theorem implies that every Lie algebra can be
embedded into $A^{(-)}$ for an appropriate associative algebra~$A$.

Leibniz algebras introduced in \cite{Loday:93,Cuvier:94} are the most popular
non-commutative analogues of Lie algebras. An algebra $(L,
[\cdot,\cdot])$ is said to be a (right) Leibniz algebra if the
product $[\cdot,\cdot]\colon L\times L\to L$ satisfies the following
(right) Leibniz identity:
\begin{equation}\label{eq:IdOfLeibnizAlgebras}
[[x,y],z]=[[x,z],y]+[x,[y,z]].
\end{equation}

To get an analogue of the Poincar\'e---Birkhoff---Witt Theorem for
Leibniz algebras, the notion of an associative dialgebra was
introduced in \cite{LodayPirashvili:93}. Namely, an associative
dialgebra is a linear space $D$ with two bilinear operations
$\vdash,\dashv\colon D\times D\to D$ satisfying certain axioms. The
new product $[x,y]=x\dashv y-y\vdash x$, % Changed, c.f. left/right Leibniz identities
$x,y\in D$, satisfies (\ref{eq:IdOfLeibnizAlgebras}), so $D$ is a
Leibniz algebra with respect to this new product. The Leibniz
algebra obtained is denoted by $D^{(-)}$. As it was shown in
\cite{Loday:01,AymonGrivel:03}, every Leibniz algebra can be
embeddable into $D^{(-)}$  for an appropriate associative
dialgebra~$D$.

Another important class of nonassociative algebras is the class of
Jordan algebras defined by the commutativity and Jordan identity
 $(x^2y)x=x^2(yx)$.
It is well-known that if $A$ is an associative algebra over a
field on characteristic  $\not=2$ then $A$ with respect to the new
product $a\circ b=\frac{1}{2}(ab+ba)$ is a Jordan algebra denoted by
$A^{(+)}$. For Jordan algebras, the analogue of the
Poincar\'e---Birkhoff---Witt theorem is not true: There exist Jordan
algebras that can not be embedded into $A^{(+)}$ for any associative
algebra~$A$. Therefore, the following notion makes sense: If a
Jordan algebra $J$ is a subalgebra of $A^{(+)}$ for some associative
algebra $A$ then it is said to be a special Jordan algebra.

The notion of a Jordan dialgebra was introduced in \cite{Kol:08} as
a particular example of a general algebraic definition of what is a
variety of dialgebras. This general operadic approach leads to three
identities defining the variety of Jordan dialgebras. Independently,
the notion of quasi-Jordan algebra emerged in
\cite{VelasquezFelipe:08} as the variety of some non-commutative
analogues of Jordan algebras. Namely, if one considers an
associative dialgebra $D$ with respect to a new product $x\circ
y=\frac{1}{2}(x\dashv y+y\vdash x)$, $x,y\in D$, then the algebra
obtained is a quasi-Jordan algebra. In \cite{VelasquezFelipe:08},
two identities were stated to define the variety of quasi-Jordan
algebras. Later in \cite{Br:08}, the third (missing) one was
noticed, so the notions of quasi-Jordan algebras and Jordan
dialgebras went to coincidence.

In \cite{Br:09}, the natural notions of a special Jordan dialgebra
and of a special identity (s-identity, for short) were introduced. An
s-identity of Jordan dialgebras is an identity which holds in all
special Jordan dialgebras but does not hold in some Jordan
dialgebra. In this note, we show the correspondence between multilinear
s-identities of Jordan algebras and Jordan dialgebras (Theorem
\ref{thm:CorrespSId}). In particular, one of the main results of
\cite{Br:09} follows from this theorem.

Also, several natural problems were posed in \cite{Br:09}: How to
generalize the classical statements known for Jordan algebras to the
case of dialgebras. This paper is devoted to the solution of all
these problems. We prove the analogues of the following theorems:
\begin{itemize}
 \item Cohn's Theorem \cite{Cohn:54}
  on the characterization of elements of free
  special Jordan algebra as symmetric elements
  of free associative algebra.
 \item Cohn's example \cite{Cohn:54} of an exceptional
  Jordan algebra which is a homomorphic image of two-generated
  special Jordan algebra. In particular, the class of special Jordan
  algebras is not a variety.
 \item Shirshov's Theorem \cite{Zhevl:78} on the speciality of
  two-generated Jordan algebra.
 \item Macdonald's Theorem \cite{Zhevl:78}
  on special identities in three variables.
\end{itemize}

The main method of study is the following. Given a Jordan dialgebra
$J$, we build two Jordan algebras $\bar J$ and $\widehat J$ as
described in \cite{Pozh:09}. The classical theorems hold for these
Jordan algebras, and their properties allow to make conclusions
about the dialgebra $J$. Moreover, the theory of conformal algebras
\cite{Kac:96} is deeply involved into considerations.

\section{PRELIMINARIES}
\subsection{Dialgebras}\label{subsec:DefDialg}
A linear space $D$ with two bilinear operations
$\vdash,\dashv\colon D\times D\to D$ is called a \emph{dialgebra}.
The base field is denoted by $\Bbbk$. A dialgebra is
\emph{associative} if it satisfies the identities
\begin{equation}\label{eq:0-DialgebraDef}
(x\mbdv y)\mbvd z=(x\mbvd y)\mbvd z,\quad x\mbdv(y\mbvd
z)=x\mbdv(y\mbdv z)
\end{equation}
and
\begin{equation}
\begin{gathered}(x,y,z)_\vdash:=(x\mbvd y)\mbvd z-x\mbvd(y\mbvd
z)=0,\\
(x,y,z)_\dashv:=(x\mbdv y)\mbdv z-x\mbdv(y\mbdv z)=0,\\
(x,y,z)_\times:=(x\mbvd y)\mbdv z-x\mbvd(y\mbdv z)=0.
\end{gathered}
\end{equation}

This class of dialgebras is well investigated in \cite{Loday:01}. Recently,
some interesting structural results on associative dialgebras were presented
in \cite{Gonzalez:10}.

A dialgebra that satisfies the identities (\ref{eq:0-DialgebraDef}), is
called a \emph{0-dialgebra.} If $D$ is a 0-dialgebra then the subspace
$D_0=\Span\{a\mbvd b-a\mbdv b\mid a,b\in D\}$ is an ideal of $D$ and
the quotient dialgebra $\bar D=D/D_0$ can be identified with an
ordinary algebra. The space $D$ may be endowed with left and right
actions of $\bar D$:
$$\bar a\cdot x=a\mbvd x,\quad x\cdot\bar a=x\mbdv a,\quad x,\,a\in
D,$$ where $\bar a$ denotes the image of $a$ in $\bar D$.

Let $A$ be an algebra that acts on a linear space $M$ via some operations
$\circ\colon A\times M\to M$ and $\circ\colon M\times A\to M$. In
this case, we can define the algebra $(A\oplus M,\circ)$, where the
product $\circ$ is given by the formula $(a+m)\circ(b+n)=ab+(a\circ
n+m\circ b)$, that is, $M\circ M=0$. The algebra obtained is called
the split null extension of $A$ by means of $M$.

We have seen before that we can define actions of the algebra $\bar
D$ on the dialgebra $D$, so the split null extension $\bar D\oplus
D$ is defined. We will denote it by $\widehat D$.

In any dialgebra $D$ a \emph{dimonomial} is an expression of the
form $w=(a_1\ldots a_n)$, where $a_1,\ldots,a_n\in D$ and
parentheses indicate some placement of parentheses with some choice
of operations. By induction we can define the \emph{central letter}
$c(w)$ of a dimonomial: if $w\in D$, then $c(w)=w$, otherwise
$c(w_1\mbvd w_2)=c(w_2)$ and $c(w_1\mbdv w_2)=c(w_1)$. Let
$c(w)=a_k$. If $D$ is 0-dialgebra, then $w=(a_1\mbvd\ldots\mbvd
a_{k-1}\mbvd a_k\mbdv a_{k+1}\mbdv\ldots\mbdv a_n)$ for the same
parenthesizing as in $(a_1\ldots a_n)$. We will denote this $w$ by
$(a_1\ldots a_{k-1}\dot a_k a_{k+1}\ldots a_n)$. In an associative
dialgebra parenthesizing does not matter, so it is reasonable to use
the notation $w=a_1\ldots a_{k-1}\dot a_k a_{k+1}\ldots a_n$, where
the dot indicates the central letter.

Let $X$ be a set of generators. Obviously, the basis of the
free dialgebra $\Di\Alg\,\langle X\rangle$ generated by $X$ consists
of dimonomials with a free placement of parentheses and a free
choice of operations. It is clear that the basis of the free
0-dialgebra $\Di\Alg 0\, \langle X\rangle$ is the set of dimonomials
$(a_1\ldots a_{k-1}\dot a_k a_{k+1}\ldots a_n)$ where $k=1,\ldots,n$
and $a_1,\ldots,a_n\in X$. Finally, the basis of the free
associative dialgebra $\Di\As\,\langle X\rangle$ consists of
dimonomials $a_1\ldots a_{k-1}\dot a_k a_{k+1}\ldots a_n$ (see
\cite{Loday:01}).

If $X=\{x_1,\ldots,x_n\}$ then every dipolynomial $f\in
\Di\As\,\langle X\rangle$ can be presented as a sum
$f=f_1+\ldots+f_n$, where each $f_i$ collects all those dimonomials
with central letter $x_i$, $i=1,\dots, n$.

\subsection{Jordan dialgebras}
Let us consider the class of Jordan dialgebras over a field $\Bbbk$
such that $\charact\Bbbk\not=2,3$. In this case, the variety of
Jordan algebras $\Jord$ over the field $\Bbbk$ is defined by the
following multilinear identities $$x_1x_2=x_2x_1,\
J(x_1,x_2,x_3,x_4)=0,$$ where
\begin{gather*}
J(x_1,x_2,x_3,x_4)=x_1(x_2(x_3x_4))+(x_2(x_1x_3))x_4+x_3(x_2(x_1x_4))\\
-(x_1x_2)(x_3x_4)-(x_1x_3)(x_2x_4)-(x_3x_2)(x_1x_4)
\end{gather*}
is the Jordan identity in a multilinear form \cite{Zhevl:78}.

Hence using the general definition of a variety of dialgebras
\cite{Kol:08} we obtain that the class of Jordan dialgebras is
defined by two 0-identities (\ref{eq:0-DialgebraDef}) and the following
identities
\begin{equation}\label{eq:IdOfJordanDialgebras}
\begin{gathered}
x_1\mbvd x_2=x_2\mbdv x_1, \\
J(\dot{x}_1,x_2,x_3,x_4)=0,\quad J(x_1,\dot{x}_2,x_3,x_4)=0, \\
J(x_1,x_2,\dot{x}_3,x_4)=0,\quad J(x_1,x_2,x_3,\dot{x}_4)=0.
\end{gathered}
\end{equation}

The variety of Jordan dialgebras is denoted $\Di\Jord$. We can
express both operations in a Jordan dialgebra through one operation:
$a\mbvd b=ab$, $a\mbdv b=ba$. Then an ordinary algebra arises that
is a noncommutative analogue of a Jordan algebra. The corresponding
variety is defined by the system of identities
$$[x_1 x_2]x_3= 0, \quad (x_1^2,x_2,x_3)=2(x_1,x_2,x_1x_3), \quad
x_1(x_1^2 x_2)=x_1^2(x_1 x_2),$$ that is equivalent to identities
(\ref{eq:IdOfJordanDialgebras}).

Such algebras are investigated in \cite{Br:08,Br:09,GubKol:09}.

\subsection{Conformal algebras}
The notion of a conformal algebra over a field of zero
characteristic was introduced by V.~G.~Kac \cite{Kac:96} as a tool
of the conformal field theory in mathematical physics. Over a field
of an arbitrary characteristic, it is reasonable to use the
following equivalent definition \cite{Kol:08}: a \emph{conformal
algebra} is a linear space $C$ endowed with a linear mapping
$T\colon C\to C$ and a set of bilinear operations ($n$-products)
$(\cdot\co{n}\cdot)\colon C\times C\to C$. For all $a,b\in C$ there
exist just a finite number of elements $n\in \mathbb{Z}^+$ such that
$a\co{n} b\not=0$ (locality property). In addition, these operations
satisfy the following properties:
\begin{gather*}
Ta\co{n}b=a\co{n-1}b,\ n\ge 1,\quad Ta\co{0}b=0,\\
T(a\co{n}b)=a\co{n}Tb+Ta\co{n}b,\ n\ge 0,
\end{gather*}
for all $a,b\in C$.

Let $\Var$ be a variety of ordinary algebras. It was defined in
\cite{Roitman:99} what is the corresponding variety of conformal
algebras when $\charact\Bbbk = 0$. Namely, given a conformal algebra $C$
one may consider $\Coeff C =\Bbbk[t,t^{-1}]\otimes _{\Bbbk[T]} C$,
where $\Bbbk[t,t^{-1}]$ is a right $\Bbbk[T]$-module defined by
$f(t)\,T = -f'(t)$, $f(t)\in \Bbbk[t,t^{-1}]$. Denote elements of
$\Coeff C$ by $a(n):=t^n\otimes _{\Bbbk[T]}a$, where $n\in \mathbb{Z}$, $a\in C$.
The multiplication on $\Coeff C$ is given by the formula
$$
a(n)\,b(m)=\sum_{s\ge 0} (-1)^s (n+m-s) \frac{n!}{(n-s)!} a\co{s}b.
$$
In \cite{Roitman:99} the definition was given: $C$ is conformal algebra
corresponding to a variety $\Var$ ($\Var$-conformal algebra) iff $\Coeff C \in \Var$.
In \cite{Kol:06} the notion of $\Var$-conformal algebra was rephrased in terms of
pseudo-algebras, that works for nonzero characteristic of $\Bbbk$.
Since we use the term "conformal algebra" for a pseudo-algebra over
$\Bbbk[T]$ in this paper, it is possible to define the class of
these objects corresponding to the variety $\Var$ of ordinary
algebras. The class of all $\Var$-conformal algebras is closed under subalgebras
and homomorphic images, but it is not closed under (infinite) direct products.
Therefore, this is not a real variety of algebraic system. We will denote it by $\Var_\mathrm{C}$.

It was also observed in \cite{Kol:08} that if $C \in \Var_\mathrm{C}$, then the space $C$ can be endowed with
a structure of a dialgebra by means of
$$a\mbvd b=a\co{0}b,\quad a\mbdv b=\sum_{s\ge 0} T^s(a\co{s}b).$$
The dialgebra obtained is denoted by $C^{(0)}$, it belongs to the
variety $\Di\Var$.

The simplest example of a conformal algebra can be constructed as
follows. Let $A$ be an ordinary algebra, then a conformal product is uniquely
defined on $\Bbbk[T]\otimes A$ by the following formulas for $a,b\in
A$:
$$a\co{n}b=\begin{cases}ab, & n=0,\\ 0, & n>0.\end{cases}$$

The conformal algebra obtained is denoted by $\Cur A$ and is called
a current conformal algebra. If an algebra $A \in \Var$, then $\Cur A \in \Var_\mathrm{C}$. In
the language of category theory, we can say that $\Cur$ is a functor
from the category of algebras to the category of conformal algebras.
If $\phi\colon A\to B$ is a homomorphism of algebras, then the
mapping $\Cur\phi\colon\Cur A\to\Cur B$ acting by the rule
$\Cur\phi(f(T)\otimes a) = f(T)\otimes \phi(a)$ is a morphism
of conformal algebras.

In \cite{GubKol:09} it was proved that an arbitrary dialgebra $D$ is
embedded into the dialgebra $(\Cur\widehat D)^{(0)}$.

\subsection{Notation for varieties of algebras and dialgebras}
An arbitrary variety of ordinary algebras we denote $\Var$, the free
algebra in this variety generated by a set $X$ is denoted by
$\Var\,\langle X\rangle$. The corresponding variety of dialgebras is
denoted by $\Di\Var$, the free dialgebra is denoted by
$\Di\Var\,\langle X\rangle$. The denotation for concrete varieties
is analogous, for example, $\Jord$ is the variety of Jordan algebras,
$\Di\Jord\,\langle X\rangle$ is the free Jordan dialgebra.

\section{SPECIAL JORDAN DIALGEBRAS}
In this section $\charact\Bbbk\not=2$. This is necessary to define
the Jordan product correctly.
\subsection{Special and exceptional Jordan dialgebras}
Let $D$ be an associative dialgebra. If we define on the set $D$ new
operations
\begin{equation}\label{eq:QuasiJordanProduct} a\mbvdjord
b=\frac{1}{2}(a\mbvd b+b\mbdv a),\ a\mbdvjord b=\frac{1}{2}(a\mbdv
b+b\mbvd a)\end{equation} then we obtain a new dialgebra which is denoted
by $D^{(+)}$. It is easy to check that this dialgebra is Jordan
\cite{Br:08}.

A dialgebra $J$ is called \emph{special}, if $J\hookrightarrow
D^{(+)}$ for some associative dialgebra $D$. Jordan dialgebras that are not
special we call \emph{exceptional}. Further, we will denote the
operations in a special Jordan dialgebra through $\mbvdjord$ and
$\mbdvjord$. These operations are expressed through associative
operations by the formula (\ref{eq:QuasiJordanProduct}).

The definition of special Jordan dialgebras has been introduced by
the analogy with ordinary algebras, where a Jordan algebra $J$ is
called special, if $J\hookrightarrow A^{(+)}$ for some associative
algebra $A$ and the product in $A^{(+)}$ is given by the formula
\begin{equation}\label{eq:JordanProduct}
a\circ b=\frac{1}{2}(ab+ba).
\end{equation}

Let now $D$ be an associative dialgebra. The mapping $*\colon D\to
D$ is called an \emph{involution} (involution of the second type
\cite{Pozh:09}) of the dialgebra $D$, if $*$ is linear and
\begin{equation}\label{eq:DefOfInvolution}
(a^*)^*=a,\quad (a\mbvd b)^*=b^*\mbdv a^*,\quad (a\mbdv
b)^*=b^*\mbvd a^*
\end{equation}
for all $a$, $b\in D$.

The set $H(D,*)=\{x\in D\mid x=x^*\}$ of symmetric elements with respect
to $*$ is closed relative to operations
(\ref{eq:QuasiJordanProduct}). This set is a subalgebra of the algebra
$D^{(+)}$. So, $H(D,*)$ is a special Jordan dialgebra.

We now construct an example of an exceptional Jordan dialgebra.
\begin{prop}\label{prop:ExampleExceptionalDialgebra} Let $(J,\circ)$
be an exceptional Jordan algebra and suppose  the condition $x\circ J=0$,
$x\in J$, implies $x=0$. Then $J$ as a dialgebra with
equal operations $x\mbvdjord y:=x\circ y$ and $x\mbdvjord y:=x\circ
y$ is an exceptional Jordan dialgebra.
\end{prop}
\begin{proof}
Assume the opposite. Let $J\hookrightarrow D^{(+)}$ where
$(D,\vdash,\dashv)$ is an associative dialgebra and the product in
$D^{(+)}$ is given by the formula (\ref{eq:QuasiJordanProduct}).
Consider $I=\Span\{a\mbvd b-a\mbdv b\mid a,\,b\in D\}$ that is an
ideal of $D$. Then $\bar D=D/I$ is an ordinary associative algebra
and $\phi\colon D^{(+)}\to\bar D^{(+)}$ is the natural
homomorphism of a Jordan dialgebra on its quotient algebra. The
composition of the embedding $\hookrightarrow$ and $\phi$ is a
homomorphism too, we denote this homomorphism through $\psi$. It is
clear that $K:=\ker\psi$ is an ideal of $J$. Since $\psi$ is a
restriction $\phi$ on $J$ so $K=\ker\psi\subseteq\ker\phi=I$.
We have $I\mbvd J=J\mbdv I=0$, this is a consequence of the
0-identity. Hence $I\circ J=I\mbvdjord J=\frac{1}{2}(I\mbvd J+J\mbdv
I)=0$, from conditions of the proposition we obtain $I=0$ therefore
and $K=0$. So $\psi$ is an embedding and $J\hookrightarrow\bar
D^{(+)}$, i. e., $J$ is exceptional.
\end{proof}

Let $\mathbf{C}$ be the Cayley-Dickson algebra over the field
$\Bbbk$, $\charact\Bbbk\not=2$. Consider an algebra $H(\mathbf{C}_3)$
of those $3\times 3$ matrices over $\mathbf{C}$ that are symmetric relative
the involution in $\mathbf{C}$. This is so called Albert algebra.
It is well-known that $J=H(\mathbf{C}_3)$ is a simple exceptional
Jordan algebra, so $J$ satisfies
the conditions of Proposition \ref{prop:ExampleExceptionalDialgebra}. Therefore,

\begin{cor}
The Albert algebra is exceptional as a Jordan dialgebra.
\end{cor}

\subsection{Symmetric and Jordan polynomials}
Let $\Alg\,\langle X\rangle$ be a free non-associative algebra
generated by $X$, $\As\,\langle X\rangle$ be a free associative
algebra, $\Di\Alg\,\langle X\rangle$ be a free non-associative
dialgebra, $\Di\As\,\langle X\rangle$ be a free associative dialgebra
\cite{Loday:01}. Products in $\Alg\,\langle X\rangle$ and
$\As\,\langle X\rangle$, also in $\Di\Alg\,\langle X\rangle$ and
$\Di\As\,\langle X\rangle$ are denoted identically. There is no
confusion because by the origin of elements it is clear which the
product we mean. Fix $z\in X$ and introduce the following mappings.

A mapping $\mathcal{J}\colon\Alg\,\langle X\rangle\to \As\,\langle
X\rangle$ is defined by linearity, on non-associative words it is defined by
induction on a length: if $x\in X$ then
$\mathcal{J}(x)=x$; if $uv\in \Alg\,\langle X\rangle$ then
$\mathcal{J}(uv)=\frac{1}{2}(\mathcal{J}(u)\mathcal{J}(v)+\mathcal{J}(v)\mathcal{J}(u))$.
So, the value of $\mathcal{J}$ on a non-associative polynomial $f$ is equal to an
associative polynomial obtained from $f$ by means of rewriting
all products in $f$ as Jordan ones by the formula (\ref{eq:JordanProduct}).
By analogy, in the case of
dialgebras a mapping $\mathcal{J}_{\Di}\colon\Di\Alg\,\langle
X\rangle\to \Di\As\,\langle X\rangle$ is defined. It is linear, it
acts identically on $x\in X$ and

$$
\begin{gathered}
\mathcal{J}_{\Di}(u\mbvd
v)=\frac{1}{2}(\mathcal{J}_{\Di}(u)\mbvd
\mathcal{J}_{\Di}(v)+\mathcal{J}_{\Di}(v)\mbdv
\mathcal{J}_{\Di}(u)),\\
\mathcal{J}_{\Di}(u\mbdv
v)=\frac{1}{2}(\mathcal{J}_{\Di}(u)\mbdv
\mathcal{J}_{\Di}(v)+\mathcal{J}_{\Di}(v)\mbvd
\mathcal{J}_{\Di}(u)).
\end{gathered}
$$

Introduce the following notation $$\Alg_z\,\langle X\rangle=\{\Phi\in
\Alg\,\langle X\rangle \mid \Phi=\sum f_i,\ f_i\text{~--- monomials,
}\deg_z f_i = 1\},$$
$$\Di\Alg_z\,\langle X\rangle=\{\Phi\in
\Di\Alg\,\langle X\rangle \mid \Phi=\sum f_i,\ f_i\text{~---
dimonomials, }\deg_z f_i = 1,\ c(f_i)=z\},$$ where $c(f_i)$
stands for the central letter of a dimonomial $f_i$. A mapping
$\Psi^z_{\Alg}\colon\Alg_z\,\langle X\rangle\to \Di\Alg_z\,\langle
X\rangle$ places signs of dialgebraic operations in a
non-associative polynomial in such a way that every sign points to $z$. By induction it can be
defined as follows: $\Psi^z_{\Alg}(z)=z$; if $z$ is contained
by $u$ then $\Psi^z_{\Alg}(uv)=\Psi^z_{\Alg}(u)\mbdv v^\dashv$; if
$z$ is contained by $v$ then $\Psi^z_{\Alg}(uv)=u^\vdash\mbvd
\Psi^z_{\Alg}(v)$. There we introduce two mappings
${}^\vdash,{}^\dashv\colon\Alg\,\langle X\rangle\to\Di\Alg\,\langle
X\rangle$. The first mapping maps a word $u$ to $u^\vdash$ where
the word $u^\vdash$ has the same multipliers as $u$ and
all signs of operations point to the right. In $v^\dashv$ all signs
of operations point to the left respectively. Further in Lemmas
\ref{lemma2} and \ref{lemma:CommutOperJPhi} we use mappings
${}^\vdash,{}^\dashv\colon\As\,\langle X\rangle\to\Di\As\,\langle
X\rangle$ which are defined and denoted in a similar way.

Analogously, we may define the sets $\As_z\,\langle X\rangle$, $\Di\As_z\,\langle
X\rangle$ and a mapping $\Psi^z_{\As}\colon\As_z\,\langle
X\rangle\to \Di\As_z\,\langle X\rangle$.

Define the following sets:
$$\SJ\,\langle X\rangle=\mathcal{J}(\Alg\,\langle X\rangle),$$
$$\Di\SJ\,\langle X\rangle=\mathcal{J}_{\Di}(\Di\Alg\,\langle X\rangle).$$

It is well-known that $\SJ\,\langle X\rangle$ is the free special Jordan dialgebra.
In fact, $\Di\SJ\,\langle X\rangle$ is the free special Jordan dialgebra (see
Lemma \ref{lemma:FreeSpecJordDialgebra} below). From the definition of the mapping $\mathcal{J}$ it is clear that
$\SJ\,\langle X\rangle$ is a subalgebra in ${\As\,\langle
X\rangle}^{(+)}$ generated by the set $X$. Similarly,
$\Di\SJ\,\langle X\rangle\hookrightarrow {\Di\As\,\langle
X\rangle}^{(+)}$.

An element from $\As\,\langle X\rangle$ is called a \emph{Jordan
polynomial} if it belongs to $\SJ\,\langle X\rangle$. By analogy, an
element from $\Di\As\,\langle X\rangle$ is called a \emph{Jordan
dipolynomial} if it belongs to $\Di\SJ\,\langle X\rangle$.

\begin{lemma}\label{lemma2}
For arbitrary $u\in\Di\As\,\langle X\rangle$, $v\in\Alg\,\langle
X\rangle$ we have
$$u\mbdv \mathcal{J}(v)^\dashv=u\mbdv \mathcal{J}_{\Di}(v^\dashv)=u\mbdv \mathcal{J}_{\Di}(v^\vdash),$$
$$\mathcal{J}(v)^\vdash\mbvd u=\mathcal{J}_{\Di}(v^\vdash)\mbvd u=\mathcal{J}_{\Di}(v^\dashv)\mbvd u.$$
\end{lemma}
\begin{proof}
Use an induction on the length of the word $v$. A base is evident.
Let $v=v_1 v_2$. Then
\begin{gather*}
u\mbdv\mathcal{J}(v)^{\dashv} =u\mbdv\mathcal{J}(v_1
v_2)^{\dashv}=\frac{1}{2}u\mbdv(\mathcal{J}(v_1)^{\dashv}\mbdv\mathcal{J}(v_2)^{\dashv}
+\mathcal{J}(v_2)^{\dashv}\mbdv\mathcal{J}(v_1)^{\dashv}) \\
=\frac{1}{2}u\mbdv(\mathcal{J}_{\Di}(v_1^{\dashv})
\mbdv\mathcal{J}_{\Di}(v_2^{\dashv})
+\mathcal{J}_{\Di}(v_2^{\dashv})\mbvd\mathcal{J}_{\Di}(v_1^{\dashv}))=u\mbdv
\mathcal{J}_{\Di}(v_1^{\dashv}\mbdv
v_2^{\dashv})=u\mbdv\mathcal{J}_{\Di}(v^{\dashv}).
\end{gather*}
All remaining equalities are proved in the same way.
\end{proof}

\begin{lemma}\label{lemma:CommutOperJPhi}
For all $\Phi\in \Alg_z\,\langle X\rangle$ we have
$$\Psi^z_{\As}(\mathcal{J}(\Phi))=\mathcal{J}_{\Di}(\Psi^z_{\Alg}(\Phi)).$$
\end{lemma}
\begin{proof}
Since all mappings are linear, it is enough to prove the statement
for the case when $\Phi$ is a word. If $\Phi=z$ then the claim is
evident. If $\Phi=uv$ then $z$ can belong to either $u$ or $v$. Let
$z$ belongs to $u$. Then using Lemma \ref{lemma2} we obtain
\begin{gather*}
\Psi^z_{\As}(\mathcal{J}(\Phi))=\Psi^z_{\As}(\frac{1}{2}[\mathcal{J}(u)\mathcal{J}(v)
+\mathcal{J}(v)\mathcal{J}(u)]) \\
=\frac{1}{2}[\Psi^z_{\As}(\mathcal{J}(u))\mbdv\mathcal{J}(v)^{\dashv}
+\mathcal{J}(v)^{\vdash}\mbvd\Psi^z_{\As}(\mathcal{J}(u))]   \\
=\frac{1}{2}[\mathcal{J}_{\Di}(\Psi^z_{\Alg}(u))\mbdv\mathcal{J}_{\Di}(v^{\dashv})
+\mathcal{J}_{\Di}(v^{\dashv})\mbvd\mathcal{J}_{\Di}(\Psi^z_{\Alg}(u))] \\
=\mathcal{J}_{\Di}(\Psi^z_{\Alg}(u)\mbdv v^{\dashv})
=\mathcal{J}_{\Di}(\Psi^z_{\Alg}(\Phi)).
\end{gather*}
The case when $z$ belongs to $v$ is proved similarly.
\end{proof}

We recall about the quotient $\bar D=D/D_0$ that has been defined in Section
\ref{subsec:DefDialg}. This quotient compares every 0-dialgebra with
an ordinary algebra. The quotient $\bar D$ of a dialgebra $D$ generated by a set $X=\{x_i \mid i\in I\}$
is an algebra generated by the set $\bar X =\{\bar x_i \mid i\in I\}$.
In order to simplify notation, we will further denote $\bar x\in \bar X$ by~$x$. Following this convention we obtain, for example,
$\overline{\Di\As\,\langle X\rangle}=\As\,\langle X\rangle$.

\begin{prop}\label{prop:SJDiSJ}
Let $f\in \Di\As_z\,\langle X\rangle$. Then
$$f\in \Di\SJ\,\langle X\rangle \Leftrightarrow \bar f\in \SJ\,\langle X\rangle.$$
\end{prop}
\begin{proof}
"$\Rightarrow$". Let $f\in \Di\SJ\,\langle X\rangle$ that is
$f=\mathcal{J}_{\Di}(\Phi)$ for some $\Phi\in\Di\Alg\,\langle
X\rangle$. Then $\bar
f=\overline{\mathcal{J}_{\Di}(\Phi)}=\mathcal{J}(\bar\Phi)$, so
$\bar f\in \SJ\,\langle X\rangle$. There we have used the equality
$\overline{\mathcal{J}_{\Di}(\Phi)}=\mathcal{J}(\bar\Phi)$ which is
easy to prove by induction on the length of $\Phi$.

"$\Leftarrow$". Let $\bar f\in \SJ\,\langle X\rangle$ that is $\bar
f=\mathcal{J}(\Phi)$ for some $\Phi\in\Alg\,\langle X\rangle$.
Since the degrees on variables do not change when we apply $\mathcal{J}$, we obtain
$\Phi\in\Alg_z\,\langle X\rangle$. Thereby, $\Phi\in\Alg_z\,\langle
X\rangle$. By Lemma \ref{lemma:CommutOperJPhi} we obtain
$\mathcal{J}_{\Di}(\Psi^z_{\Alg}(\Phi))=\Psi^z_{\As}(\mathcal{J}(\Phi))
=\Psi^z_{\As}(\bar f)=f$, the last equality in the sequence is true
because $f\in\Di\As_z\,\langle X\rangle$. So, $f\in \Di\SJ\,\langle
X\rangle$.
\end{proof}

Consider the dialgebra
\[
\Lambda_X =  \Di\As\,\langle X\rangle /I,
\]
where $I$ is the ideal of $\Di\As\,\langle X\rangle$
generated by the set $\{ f_{x,y} = x\mbdv y + y\mbvd x \mid x,y \in
X\}$. This dialgebra is the analogue of the exterior algebra
(Grassmann algebra). Further we will identify the set $X$ and its
image $\bar X\subseteq\Lambda_X$. Following this
agreement we suppose that $\Lambda_X$ is generated by the set $X$.

\begin{thm}\label{thm:BasisGrassmanDialgebra}
Let $X$ be a linearly ordered set. Then the basis of the algebra
$\Lambda_X$ consists of words $\dot x_1x_2\dots x_k$, $k\ge
1$, $x_i\in X$, $x_2<x_3<\dots < x_k$.
\end{thm}
\begin{proof}
Use the theory of Gr\"{o}bner-Shirshov bases for associative
dialgebras developed in \cite{BokutChenLiu:08}. Let $S_0 =
\{f_{x,y}\mid x,y\in X \}$ be the initial set of defining relations.
Compositions of left product $z\mbdv f_{x,y}$ belong
to the ideal $I$ as well as compositions of right product
$f_{x,y}\mbvd z$, $x,y,z\in X$. The set of defining relations
obtained
\[
x\mbdv y + y\mbvd x; \quad x\mbdv y\mbdv z + x\mbdv z\mbdv y,\ y>z;
\quad x\mbvd y\mbvd z + y\mbvd x\mbvd z,\ x>y; \quad x\mbvd x\mbvd
y; \quad x\mbdv y\mbdv y
\]
is a Gr\"{o}bner-Shirshov basis. Reduced words are of the form
$$\dot x_1x_2\dots x_k,\ k\ge 1,\ x_2<x_3<\dots < x_k,$$ and
the set of all reduced words by \cite{BokutChenLiu:08} is a linear
basis of the algebra $\Lambda_X$.
\end{proof}

Define an involution $*$ on $\Di\As\,\langle X\rangle$ in the following way:
$$(x_{i_1}\ldots \dot{x}_{i_k}\ldots
x_{i_n})^*=x_{i_n}\ldots \dot{x}_{i_k}\ldots x_{i_1},$$ and extend to
dipolynomials by linearity. This mapping reverses words
and signs of dialgebraic operations. It is easy to check that the
mapping $*$ satisfies properties of an involution
(\ref{eq:DefOfInvolution}). Through $\Di\Herm\,\langle X\rangle$ we
denote the Jordan dialgebra $H(\Di\As\,\langle X\rangle,*)$ of symmetric
elements from $\Di\As\,\langle X\rangle$ with respect to $*$
with the product (\ref{eq:QuasiJordanProduct}).

Further $\{u\}$ denotes $ u+u^*$ where $u$ is a basic word from
$\Di\As\,\langle X\rangle$. Note that $\{u\}=\{u^*\}$.

An analogous involution on $\As\,\langle X\rangle$ we denote by $*$
too. It acts like as
$$(x_{i_1}\ldots x_{i_k}\ldots x_{i_n})^*=x_{i_n}\ldots x_{i_k}\ldots x_{i_1},$$
on monomials and extends to polynomials by linearity.

The next theorem is an analogue of the classical Cohn's Theorem
\cite[Theorems 4.1 and 4.2]{Cohn:54} that describes Jordan polynomials from $\le 3$
variables as symmetric elements of the free associative algebra.

\begin{thm}\label{thm:CohnForDialgebra}
For any set $X$ we have $\Di\SJ\,\langle
X\rangle\subseteq\Di\Herm\,\langle X\rangle$. If $|X|\le 2$
then there is an equality, if $|X|>2$ then there is a strict
inclusion. Also, for any $X$ we have that $\Di\Herm\,\langle X\rangle$ is generated by $X$
and dotted tetrads $\{\dot xyzt\}$, $\{\dot xxyz\}$, where $x,y,z,t\in X$ are distinct.
\end{thm}
\begin{proof}
"$\subseteq$" follows from the equality
$\mathcal{J}_{\Di}(\Phi)^*=\mathcal{J}_{\Di}(\Phi)$ which holds
for all $\Phi\in\Di\Alg\,\langle X\rangle$. As before, this equality
can be proved by induction on the length of $\Phi$ considering cases
$\Phi=u\mbvd v$ and $\Phi=u\mbdv v$.

Let $|X|=2$. In order to prove the equality, consider an arbitrary $f\in\Di\Herm\,\langle
x,y\rangle$, i.~e., $f\in\Di\As\,\langle x,y\rangle$ and $f=f^*$. We
need to show that $f\in\Di\SJ\,\langle x,y\rangle$. The dipolynomial
$f$ is equal to a sum of dimonomials $f=\sum f_i$. Further,
$f=\frac{1}{2}(f+f^*)=\frac{1}{2}\sum{(f_i+f_i^*)}$. Without loss of
generality we may assume $f=a+a^*$ where $a$ is a dimonomial.
Suppose $x$ is the central letter of $a$. So
$f$ can be written in a form $f=u\dot xv+v^*\dot xu^*$ where
$u,\,v\in\Di\As\,\langle x,y\rangle$ or equal to empty words.
Consider $g(x,y,z)=u\dot zv+v^*\dot zu^*\in\Di\As\,\langle
x,y,z\rangle$. Since $\bar{g}=\bar{g}^*$ then
$\bar{g}\in\SJ\,\langle x,y,z\rangle$ by the classical Cohn's
Theorem. In addition, $g\in\Di\As_z\,\langle x,y,z\rangle$ hence
Proposition \ref{prop:SJDiSJ} implies $g\in\Di\SJ\,\langle
x,y,z\rangle$. It means that there exists a dipolynomial $\Phi(x,y,z)$ such
that $g=\mathcal{J}_{\Di}(\Phi(x,y,z))$. Substituting $z:=x$ into the last
equality we obtain $f=\mathcal{J}_{\Di}(\Phi(x,y,x))$. Therefore,
$f\in\Di\SJ\,\langle x,y\rangle$. We have proved the equality
for $|X|=2$ and thus for $|X|=1$.

Let $|X|>2$. In order to prove the strict inclusion consider the
dotted tetrad $\{\dot xxyz\}=\dot xxyz+zyx\dot x\in\Di\Herm\,\langle X\rangle$ where
$x,y,z\in X$. There exists a homomorphism $\sigma\colon\Di\As\,\langle
x,y,z\rangle \to \Di\Lambda\,\langle x,y,z\rangle$ such that
$\sigma(x)=x$, $\sigma(y)=y$, $\sigma(z)=z$. All Jordan
dipolynomials of degree greater that 1 map to zero by this
homomorphism. Using the basis of $\Di\Lambda\,\langle x,y,z\rangle$
from Theorem \ref{thm:BasisGrassmanDialgebra} we obtain
$$\sigma(\{\dot xxyz\})=2\dot xxyz\neq 0.$$ (When we use Theorem
\ref{thm:BasisGrassmanDialgebra} we suppose that $x<y<z$.)
So, the dipolynomial $\{\dot xxyz\}$ does not belong to
$\Di\SJ\,\langle X\rangle$.

It is well-known that in Jordan algebras we can permute variables in a tetrad modulo $\SJ\,\langle X\rangle$. It follows from the fact that $xy = -yx + x\circ y$, hence $\{xyzt\} = -\{yxzt\} + \{(x\circ y) zt\}$ and $\{xyzt\} \in -\{yxzt\} + \SJ\,\langle X\rangle$. Placing a dot in the last equality we get that $\{\dot xyzt\} \in -\{y\dot xzt\} + \Di\SJ\,\langle X\rangle$. Therefore, we can permute the variables in a dotted tetrad together with a dot modulo $\Di\SJ\,\langle X\rangle$.

To find the generators of $\Di\Herm\,\langle X\rangle$, we consider $g(x_1,\ldots,x_n) \in \Di\Herm\,\langle X\rangle$. We can write $g=g_1+\ldots+g_n$, where each $f_i$ collects all dimonomials with central letter $x_i$. It suffices to find generators for $g_1$. There exists $f(s,x_1,\ldots,x_n) \in \Herm\,\langle X,s \rangle$ such that $g_1 = \Psi^s_\As (f) |_{s:=x_1}$. Theorem 4.1 \cite{Cohn:54} holds for $f$, hence $f$ is Jordan polynomial $\Phi$ from $X$, $s$ and tetrads, i.e., $f=\mathcal{J}(\Phi)$. Further, $g_1 = \Psi^s_\As (\mathcal{J}(\Phi)) |_{s:=x_1} \stackrel{\mathrm{L.4}}{=} \mathcal{J}_\Di (\Psi^s_\Alg(\Phi) |_{s:=x_1} )$. Therefore, $g_1$ is generated by $X$, by dotted tetrads with all the different variables and by dotted tetrads with two equal variables and with the dot over one of the equal variables. By permutation of variables these dotted tetrads can be reduced to $\{\dot xyzt\}$ and $\{\dot xxyz\}$.
\end{proof}

\subsection{Homomorphic images of special Jordan dialgebras}
In this section we construct the example of an exceptional
two-generated Jordan dialgebra which is a homomorphic image of a
special Jordan dialgebra.

Denote by $\widehat I$ the ideal of $\Di\As\,\langle x,y \rangle$
generated by the set $I$.

\begin{lemma}\label{lemma:CriterionOfQuotientSpeciality}
Let $I$ be an ideal of $\Di\SJ\,\langle X \rangle$. Then $\Di\SJ\,\langle X \rangle/I$ is special iff
$\widehat I \cap \Di\SJ\,\langle X\rangle \subseteq I$.
\end{lemma}
\begin{proof}
The proof of this lemma is completely analogous to the proof of Theorem 2.2 \cite{Cohn:54}.
\end{proof}

\begin{prop}\label{prop:CriterionOfQuotientSpecialityTwoGenerated}
Let $I$ be an ideal of $\Di\SJ\,\langle x,y \rangle$ is generated by elements $u_i$. Then $\Di\SJ\,\langle x,y \rangle/I$ is special iff
$\{u_i \dot xxy\} \in I$ and $\{u_i \dot yyx\} \in I$ for all $i$.
\end{prop}
\begin{proof}
By Theorem \ref{thm:CohnForDialgebra}, $\Di\SJ\,\langle x,y \rangle = \Di\Herm\,\langle x,y \rangle$. Lemma \ref{lemma:CriterionOfQuotientSpeciality} implies that $\Di\SJ\,\langle x,y \rangle/I$ is special iff
$\widehat I \cap \Di\Herm\,\langle x,y \rangle \subseteq I$.

"$\Rightarrow$". It is clear that $\{u_i \dot xxy\} \in \widehat I \cap \Di\Herm\,\langle x,y \rangle$, hence the condition of proposition is necessary.

"$\Leftarrow$". Suppose that $\{u_i \dot xxy\} \in I$ and  $\{u_i \dot yyx\} \in I$ for all $i$ and let $w \in \widehat I \cap \Di\Herm\,\langle x,y \rangle$. It is clear (as in Lemma 3.2 \cite{Cohn:54}) that $w$ can be written as a symmetric polynomial $f=f^*$ in the $u$'s and $x$, $y$ which is linear homogenious in the $u$'s. We now regard $x$, $y$, $u_i$ as independent. Because $f\in\Di\Herm\,\langle x,y,u_i \rangle$, it can by Theorem \ref{thm:CohnForDialgebra} be expressed as Jordan dipolynomial $\Phi$ in $x$, $y$, $u_i$ and dotted tetrads involving this variables. Since $f$ is linear in the $u$'s so is $\Phi$ and therefore no dotted tetrad can involve more than one $u$, but it must involve at least one. By permutation of variables any such tetrad can be reduced to the form $\{u_i \dot xxy\}$ or $\{u_i \dot yyx\}$ plus Jordan dipolynomial. By hypothesis any such tetrad belong to $I$, hence every term of $\Phi$ contains at least one factor from $I$, so $\Phi\in I$. This shows that $w=f=\Phi\in I$ and this completes the proof.

\end{proof}

\begin{thm}\label{thm:Example2GeneratedExceptionalDialgebra}
Consider the special Jordan dialgebra $\Di\SJ\,\langle x,y\rangle$,
and let $I$ be its ideal generated by the element
$k = \frac{1}{2} (\dot x x + x\dot x) - \frac{1}{2}(\dot yy+y\dot
y)$. Then the quotient dialgebra $J =\Di\SJ\, \langle x,y\rangle / I$
is exceptional.
\end{thm}
\begin{proof}
Consider $f = \{kx\dot xy\}$. By Proposition \ref{prop:CriterionOfQuotientSpecialityTwoGenerated} it suffices to show that $f\notin I$.

Assume $f\in I$. Then there exists a dipolynomial
\[
 \phi (x,y,z) \in \Di\SJ\, \langle x,y,z\rangle \subset \Di\Herm\,
 \langle x,y,z\rangle
\]
such that $\phi (x,y,k) = f$. In addition, every summand
from $\phi$ contains at most one entry of $z$.

Write
$$
 \phi(x,y,z)  = \phi_1 (x,y,z) + \phi_2(x,y,z) + \dots,
 \quad
 \deg_z \phi_n  = n.
$$

The total degree of $f$ (with respect to all variables) is equal to 5, hence
$\phi_n = 0$ when $n\ge 3$. Therefore $\phi(x,y,z)  =
\phi_1 (x,y,z) + \phi_2(x,y,z)$.

Suppose $\phi_1 := \phi_{1,0} +
\phi_{1,1}+\phi_{1,2}+\phi_{1,3},$ where $\deg_x
\phi_{1,0}=0$, $\deg_x \phi_{1,1}=1$, $\deg_x\phi_{1,2}=2$,
$\deg_x \phi_{1,3}=3$; $\phi_2 := \phi_{2,0} +
\phi_{2,1}$, where $\deg_x \phi_{2,0}=0$, $\deg_x
\phi_{2,1}=1$.

After the substitution $z=k$ all summands in $\phi_{1,1}$,
$\phi_{1,3}$ and $\phi_{2,1}$ have degree 1, 3 or 5 in $x$.
All summands from $\phi_{1,0}$, $\phi_{1,2}$ and
$\phi_{2,0}$ have degree 0, 2 or 4 in $x$. Since $f$ contains
summands of only 2-nd and 4-th degree in $x$, we have
$\phi_{1,1}+\phi_{1,3}+\phi_{2,1}=0$.

Therefore, $\phi=\phi_{1,0}+\phi_{1,2}+\phi_{2,0}$.

Since $x$ is the central letter of the dipolynomial $f$, central
letters of dimonomials from $\phi$ can be variables $x$ and $z$.
Every dipolynomial from $\Di\Herm\,\langle x,y,z \rangle$ with this
property is equal to a linear combination of the next
dipolynomials:
\[
\begin{gathered}
\{ \dot xyxz \}, \ \{ xy\dot xz \}, \ \{ xyx\dot z \},\quad
\{ y\dot xxz \}, \ \{ yx \dot xz \}, \ \{ yxx\dot z \},\\
\{ \dot xxyz \}, \ \{ x\dot xyz \}, \ \{ xxy\dot z \},\quad
\{ \dot xyzx \}, \ \{ xyz \dot x \}, \ \{ xy\dot zx \},\\
\{ yz\dot xx \}, \ \{ yzx\dot x \}, \ \{ y\dot zxx \} ,\quad
\{ y\dot xzx \}, \ \{ yxz\dot x \}, \ \{ yx\dot zx \},\\
\{\dot zyyy\},\ \{y\dot zyy\},\quad \{\dot zzy\},\ \{z\dot zy\},
\{\dot zyz\}.
\end{gathered}
\]

Consequently $\phi (x,y,z)$ has the form
\[
\begin{gathered}
 \alpha_1 \{ \dot xyxz \}
+\alpha_2 \{ y\dot xxz \} +\alpha_3 \{ \dot xxyz \} +\alpha_4 \{
\dot xyzx \} +\alpha_5 \{ yz\dot xx \}
+\alpha_6 \{ y\dot xzx \} \\
+\beta_1 \{ xy\dot xz \} +\beta_2 \{ yx \dot xz \} +\beta_3 \{ x\dot
xyz \} +\beta_4 \{ xyz \dot x \} +\beta_5 \{ yzx\dot x \}
+\beta_6 \{ yxz\dot x \}  \\
+2\gamma_1 \{ xyx\dot z \} +2\gamma_2 \{ yxx\dot z \} +2\gamma_3 \{
xxy\dot z \} +2\gamma_4 \{ xy\dot zx \} +2\gamma_5 \{ y\dot zxx \}
+2\gamma_6 \{ yx\dot zx \} \\
+2\delta_1\{\dot zyyy\}+2\delta_2\{y\dot zyy\}+2\delta_3\{\dot
zzy\}+2\delta_4\{z\dot zy\}+2\delta_5\{\dot zyz\}.
\end{gathered}
\]

Substituting $z=k$ and using the equalities
\begin{gather*}
2\dot zz=(\dot xx+x\dot x-\dot yy-y\dot y)\mathbin\dashv(xx-yy) \\
=\dot xx^3+x\dot xx^2-\dot yyx^2-y\dot yx^2-\dot xxy^2-x\dot
xy^2+\dot y y^3+y\dot yy^2, \\
2z\dot z=(xx-yy)\mathbin\vdash(\dot xx+x\dot x-\dot yy-y\dot y) \\
=x^2\dot xx+x^3\dot x-x^2\dot yy-x^2y\dot y-y^2\dot xx-y^2x\dot x
+y^2\dot yy+y^3\dot y,
\end{gather*}
we obtain $\phi(x,y,k)$ is equal to
\[
\begin{gathered}
 \alpha_1 \{ \dot xyx^3 \}
+\alpha_2 \{ y\dot x x^3 \} +\alpha_3 \{ \dot xxyx^2 \} +\alpha_4 \{
\dot xyx^3 \} +\alpha_5 \{ yx^2\dot xx \}
+\alpha_6 \{ y\dot xx^3 \}         \\
- \alpha_1 \{ \dot xyxy^2 \} - \alpha_2 \{ y\dot xxy^2 \} - \alpha_3
\{ \dot x xy^3 \} - \alpha_4 \{ \dot xy^3x \} - \alpha_5 \{ y^3\dot
xx \}
- \alpha_6 \{ y\dot xy^2x \}  \\
+\beta_1 \{ xy\dot x x^2 \} +\beta_2 \{ yx \dot x x^2 \} +\beta_3 \{
x\dot xyx^2 \} +\beta_4 \{ xyx^2 \dot x \} +\beta_5 \{ yx^3\dot x \}
+\beta_6 \{ yx^3\dot x \}       \\
- \beta_1 \{ xy\dot xy^2 \} - \beta_2 \{ yx \dot xy^2 \} - \beta_3
\{ x\dot xy^3 \} - \beta_4 \{ xy^3 \dot x \} - \beta_5 \{ y^3x\dot x
\}
- \beta_6 \{ yxy^2\dot x \}  \\
+\gamma_1 \{ xyx\dot x x \} +\gamma_2 \{ yx^2\dot xx \} +\gamma_3 \{
x^2y\dot xx \} +\gamma_4 \{ xy\dot xx^2 \} +\gamma_5 \{ y\dot xx^3
\}
+\gamma_6 \{ yx\dot xx^2 \}   \\
+\gamma_1 \{ xyx^2\dot x \} +\gamma_2 \{ yx^3\dot x \} +\gamma_3 \{
x^2yx\dot x\} +\gamma_4 \{ xyx\dot xx \} +\gamma_5 \{ yx\dot xx^2 \}
+\gamma_6 \{ yx^2\dot xx \}   \\
-\gamma_1 \{ xyx\dot yy \} -\gamma_2 \{ yx^2\dot yy \} -\gamma_3 \{
x^2y\dot yy \} -\gamma_4 \{ xy\dot yyx \} -\gamma_5 \{ y\dot yyx^2
\}
-\gamma_6 \{ yx\dot yyx \}   \\
-\gamma_1 \{ xyxy\dot y \} -\gamma_2 \{ yx^2y\dot y \} -\gamma_3 \{
x^2y^2\dot y \} -\gamma_4 \{ xy^2\dot yx \} -\gamma_5 \{ y^2\dot
yx^2\} -\gamma_6 \{ yxy\dot yx \}   \\
\end{gathered}
\]

\[
\begin{gathered}
+\delta_1\{\dot xxy^3\}+\delta_1\{x\dot xy^3\}-\delta_1\{\dot
yy^4\}-\delta_1\{y\dot yy^3\}   \\
+\delta_2\{y\dot xxy^2\}+\delta_2\{yx\dot xy^2\}-\delta_2\{y\dot y
y^3\}-\delta_2\{y^2\dot y y^2\}   \\
+\delta_3\{\dot xx^3y\}+\delta_3\{x\dot xx^2y\}-\delta_3\{\dot
yyx^2y\}
-\delta_3\{y\dot yx^2y\}   \\
-\delta_3\{\dot xxy^3\}-\delta_3\{x\dot xy^3\}+\delta_3\{\dot
yy^4\}+\delta_3\{y\dot yy^3\}   \\
+\delta_4\{x^2\dot xxy\}+\delta_4\{x^3\dot xy\}-\delta_4\{x^2\dot
yy^2\}-\delta_4\{x^2y\dot yy\}   \\
-\delta_4\{y^2\dot xxy\}-\delta_4\{y^2x\dot xy\}+\delta_4\{y^2\dot
yy^2\}+\delta_4\{y^3\dot yy\}   \\
+\delta_5\{\dot xxyx^2\}+\delta_5\{x\dot xyx^2\}-\delta_5\{\dot
yy^2x^2\}-\delta_5\{y\dot yyx^2\}   \\
-\delta_5\{\dot xxy^3\}-\delta_5\{x\dot xy^3\}+\delta_5\{\dot
yy^4\}+\delta_5\{y\dot yy^3\}.
\end{gathered}
\]

This expression must coincide with $f = \{x^3\dot x y \} -
\{y^2x\dot x y\} $. In particular, a sum of all dimonomials with the
central letter $y$ must be equal to zero:
\[
\begin{gathered}
0= \gamma_1\{y\dot yxyx\} +(\gamma_2+\delta_3)\{y\dot yx^2y\}
+(\gamma_3+\gamma_5+\delta_4+\delta_5)\{y\dot yyx^2\}  \\
+\gamma_4\{xy\dot yyx\} +\gamma_6\{xy\dot yxy\} +\gamma_1\{\dot
yyxyx\} +(\gamma_2+\delta_3)\{\dot y yx^2y\}
+(\gamma_3+\delta_5)\{\dot y y^2x^2\} \\
+\gamma_4\{x\dot yy^2x\} +(\gamma_5+\delta_4)\{x^2\dot yy^2\}
+\gamma_6 \{x\dot yyxy\} +(\delta_1-\delta_3-\delta_5)\{\dot yy^4\} \\
+(\delta_1+\delta_2-\delta_3-\delta_4-\delta_5)\{y\dot yy^3\}
+(\delta_2-\delta_4)\{y^2\dot yy^2\}.
\end{gathered}
\]
All coefficients in this sum have to be zero. Solving the
obtained system we have $\gamma_2= -\delta_3$, $\gamma_3=
-\delta_5$, $\gamma_5= -\delta_4$, $\delta_1= \delta_3+\delta_5$,
$\delta_2=\delta_4$, $\gamma_1= \gamma_4= \gamma _6= 0$.

Substitute the obtained relations to $\phi(x,y,k)$ we get that
all summands with coefficients $\gamma$ and $\delta$ are eliminated.

Further, consider the remaining summands (we divide them into two groups
by $\deg_y$):
\[
\begin{gathered}
(\alpha_1 + \alpha_4) \{ \dot xyx^3 \} +(\alpha_2 + \alpha_6) \{
y\dot x x^3 \} +\alpha_3 \{ \dot xxyx^2 \}
+\alpha_5 \{ yx^2\dot xx \}    \\
+\beta_1 \{ xy\dot x x^2 \} +\beta_2 \{ yx \dot x x^2 \} +\beta_3 \{
x\dot xyx^2 \} +\beta_4 \{ xyx^2 \dot x \}
+(\beta_5 +\beta_6)\{ yx^3\dot x \}   \\
= \{x^3\dot x y\},
\end{gathered}
\]
\[
\begin{gathered}
\alpha_1 \{ \dot xyxy^2 \}+\alpha_2 \{ y\dot xxy^2 \}
+\alpha_3\{\dot x xy^3 \} + \alpha_5 \{ y^3\dot xx \}
+ \alpha_6 \{ y\dot xy^2x \}  \\
+\beta_1 \{ xy\dot xy^2 \}+\beta_2 \{ yx \dot xy^2 \}+
\beta_3\{x\dot xy^3 \}+ (\alpha_4 +\beta_4) \{ \dot xy^3x \}
+\beta_5 \{ y^3x\dot x \}+\beta_6 \{ yxy^2\dot x \} \\
= \{ y^2x\dot x y\}.
\end{gathered}
\]

The last two equalities imply $\alpha _2 =1$ and
other coefficients are equal to zero. Therefore,
\begin{gather*}
 \phi (x,y,z) = \{y\dot xxz\}-2\delta_3\{yxx\dot z\}
 -2\delta_5\{xxy\dot z\}-2\delta_4\{y\dot zxx\}  \\
 +2(\delta_3+\delta_5)\{\dot zyyy\}+2\delta_4\{y\dot zyy\}
 +2\delta_3\{\dot zzy\}+2\delta_4\{z\dot zy\}+2\delta_5\{\dot zyz\}.
\end{gather*}

By assumption this dipolynomial is Jordan. When we expand Jordan
products then the central letter is preserved, hence the
dipolynomials consisting of dimonomials from $\phi (x,y,z)$ with
the fixed central letter must be Jordan. In particular, if we
choose the central letter $x$ then the dipolynomial
$\{y\dot xxz\}$ must be Jordan, but this is not true by the proof of
Theorem \ref{thm:CohnForDialgebra}.

The contradiction obtained proves that $f\notin I$.
\end{proof}

\section{S-IDENTITIES}
In this section $\charact\Bbbk=0$, so we can perform the process of
full linearization of identities and varieties of algebras are
always defined by multilinear identities.
\subsection{Equality of varieties $\Hom\Di\SJ$ and $\Di\Hom\SJ$}
Consider a class of special Jordan dialgebras $\SJ$. The class $\SJ$
is not a variety because it is not close relative the taking of
homomorphic images. Consider the operator $\Hom$ acting on classes of
algebraic systems
$$\Hom(K)=\{A\mid A=\phi(B)\text{ for }B\in K,\phi\colon B\to A
\text{ is an epimorphism}\}.$$

It is well-known that $\Hom(\SJ)$ is a variety of algebras which we
denote $\Hom\SJ$.

We recall (see Section \ref{subsec:DefDialg}) that if
$D\in\Di\Alg 0$ then $D$ can be endowed with left and right actions
of the algebra $\bar D$ by the rules $\bar xy=x\mbvd y$, $y\bar x=y\mbdv
x$, where $x$, $y\in D$. Let $\Var$ be a variety of ordinary
algebras. In the paper \cite{Pozh:09} it is shown that $D\in\Di\Var$
if and only if $\bar D\in\Var$ and $D$ is a $\Var$-bimodule over $\bar D$
in the sense of Eilenberg, i.~e., the split null extension
$\widehat D=\bar D\oplus D$ belongs to the variety $\Var$.

In this way we can define a variety of dialgebras $\Di\Hom\SJ$ by a
variety $\Hom\SJ$.

Let $\Di\SJ$ be the class of special Jordan dialgebras. Consider the
closure $\Hom(\Di\SJ)$ of this class relative to the operator $\Hom$.
The variety obtained we denote by $\Hom\Di\SJ$.

The purpose of this section is to show that
$\Hom\Di\SJ=\Di\Hom\SJ$.

\begin{lemma}\label{lemma:FreeSpecJordDialgebra}
$\Di\SJ\,\langle X\rangle$ is a free algebra in the variety
$\Hom\Di\SJ$.
\end{lemma}
\begin{proof}
Let $J'\in\Hom\Di\SJ$ be a homomorphic image of $J\in\Di\SJ$,
$D$ be an associative dialgebra such that $J\hookrightarrow D^{(+)}$.
We have the following commutative diagram
$$\begin{CD}
J'   @<\text{на}<< J         @>\subseteq>> D     \\
@AAA          @AAA             @AAA      \\
X    @>\subseteq>> \Di\SJ\langle X\rangle @>\subseteq>>
\Di\As\langle X\rangle
\end{CD}$$

We have $X\subseteq J'$. Consider some preimages of elements of the set $X$
with respect to the mapping $J\to J'$. Since $J\subseteq D$, we obtain the embedding of $X$ into $D$.
By the universal property of $\Di\As\langle X\rangle$ there exists an unique
homomorphism $\Di\As\langle X\rangle\to D$ such that its
restriction to $\Di\SJ\langle X\rangle$ is the homomorphism
$\Di\SJ\langle X\rangle\to J$. The last homomorphism in a
composition with the mapping $J\to J'$ gives the required
homomorphism $\Di\SJ\langle X\rangle\to J'$.
\end{proof}

A \emph{bar-unit} of a 0-dialgebra $D$ is an element $e\in D$ such
that $x\mbdv e=e\mbvd x=x$ for every $x\in D$ and $e$
belongs to the associative center of $D$ that is
$$(x,e,y)_\times=(e,x,y)_\dashv=(x,y,e)_\vdash=0$$ for all $x$,
$y\in D$.

\begin{prop}[{Pozhidaev \cite[Theorem 2.2]{Pozh:09}}]\label{prop:EmbWithBarUnit}
For every associative dialgebra $D$ there exists an associative
dialgebra $D_e$ with the bar-unit $e$ such that $D\hookrightarrow D_e$.
\end{prop}

\begin{lemma}\label{lemma:SpecUnitEmb}
Let $J$ be a special Jordan dialgebra. Then there exists a special Jordan
dialgebra $J_e$ such that  $J\hookrightarrow J_e$ and $\bar e$
is a unit in the algebra $\bar J_e$.
\end{lemma}
\begin{proof}
By the defintion of a special Jordan dialgebra it follows that
$J=(J,\vdjord,\dvjord)$ is embedded into $D^{(+)}$ for
some associative dialgebra $D=(D,\vdash,\dashv)$. By Proposition
\ref{prop:EmbWithBarUnit} we have an embedding
$D^{(+)}\hookrightarrow D_e^{(+)}$ where $e$ is a bar-unit in $D_e$.
Therefore, $J_e=D_e^{(+)}$ is the required dialgebra.

Further, $e\mbvd x=x\mbdv e=x$ holds for every $x\in D_e$, so in $J_e$ we
have $e\mbvdjord x=\frac{1}{2}(e\mbvd x+x\mbdv e)=x$, $x\mbdvjord
e=x$. Hence $\bar e\bar x=\bar x\bar e=\bar x$ in the quotient algebra
$\bar J_e$, so $\bar e$ is a unit in $\bar J_e$.
\end{proof}

\begin{lemma}\label{lemma:SpecFact}
Let $J$ be a special Jordan dialgebra and such that $\bar{J}$ contains a unit.
Then $\bar{J}$ is special.
\end{lemma}
\begin{proof}
Let $D$ be an associative dialgebra such that $J\hookrightarrow
D^{(+)}$. Denote $\langle D,D\rangle:=\Span\{a\mbvd b-a\mbdv b\mid
a,\,b\in D\}$, $[J,J]:=\Span\{a\mbvdjord b-a\mbdvjord b\mid a,\,b\in
J\}$. As before $\bar D=D/\langle D,D\rangle$ is an associative
algebra. Since $J\subseteq D$ we have $[J,J]\subseteq\langle
D,D\rangle$. Then the homomorphism $\phi\colon \bar J\to\bar
D^{(+)}$ is well-defined by the rule $x+[J,J]\mapsto x+\langle D,D\rangle$.

$$\begin{CD}
J @>\subseteq>> D \\
@VVV  @VVV \\
\bar J @>\phi>> \bar D
\end{CD}$$

It is evident that $\phi$ is injective if and only if $\langle
D,D\rangle\cap J=[J,J]$.

Let $x\in\langle D,D\rangle\cap J$. Then $x\mbvd y=y\mbdv x=0$ for
every $y\in D$, hence $x\mbvdjord y=\frac{1}{2}(x\mbvd y+y\mbdv
x)=0$ in $J$ and $\bar x\bar y=\bar 0$ in $\bar J$. Take $\bar
y=1\in\bar J$ and obtain $\bar x=\bar 0$, i.~e., $x\in [J,J]$. So,
$\phi$ is injective and $\bar J$ is special.
\end{proof}

Let $J$ be a Jordan algebra, $A$ be an associative algebra with a
unit, then a homomorphism from $J$ to $A^{(+)}$ is called an
\emph{associative specialization} $\sigma\colon J\to A$. This is
a linear mapping such that
$$\sigma(ab)=\frac{1}{2}(\sigma(a)\sigma(b)+\sigma(b)\sigma(a))$$
for all $a,b\in J$.

Two associative specializations are called \emph{commuting} if
$[\sigma_1(a),\sigma_2(b)]=0$ for all $a,b\in J$.

A bimodule $M$ over $J$ is \emph{special} if there exists an
embedding of $M$ into a bimodule $N$ such that if $v\in N$, $a\in J$
then
\begin{equation}\label{eq:SpecBiModCond}
a\cdot v=\frac{1}{2}(\sigma_1(a)v+\sigma_2(a)v),
\end{equation}
where $\sigma_1$, $\sigma_2$ are commuting associative
specializations of $J$ into $\mathrm{Hom}(N,N)$.

\begin{thm}[{Jacobson \cite[theorem II.17]{Jacobson:68}}]\label{thm:SpecSplitNullExtCrit}
Let $J$ be a special Jordan algebra, $M$ be a bimodule over $J$.
Then the bimodule $M$ is special if and only if the split null
extension $J\oplus M$ is a special Jordan algebra.
\end{thm}

\begin{lemma}\label{lemma:SpecSplitNullExt}
Let $J$ be a special Jordan dialgebra and $\bar{J}$ be a special
Jordan algebra. Then $\widehat{J}=\bar{J}\oplus J$ is special too.
\end{lemma}
\begin{proof}
Since $J=(J,\vdjord,\dvjord)$ is special, we have
$J\hookrightarrow D^{(+)}$ where $D={(D,\vdash,\dashv)}$ is an
associative dialgebra. The dialgebra $J$ is a $\bar
J$-bimodule: $\bar a\cdot v=a\mbvdjord v =v\mbdvjord a=v\cdot\bar
a$, where $\bar a\in\bar J$, $v\in J$.

We prove that the bimodule $J$ over the special Jordan algebra $\bar J$
is special. The bimodule $J$ is embedded into $D$ and $D$ is a $\bar
J$-bimodule too. Consider mappings $\sigma_1,\,\sigma_2\colon\bar
J\to\mathrm{Hom}(D,D)$ defined by the rule
$$\sigma_1(\bar a)\colon d\mapsto a\mbvd d\in D,\,
\sigma_2(\bar a)\colon d\mapsto d\mbdv a\in D,\quad d\in D,\,a\in
J\subseteq D.$$

These mappings are well-defined. We show that they are associative
specializations. Indeed for every $\bar{a\vphantom b},\,\bar
b\in\bar J$, $d\in D$
\begin{gather*}
\sigma_1(\bar{a\vphantom b}\bar b)(d) =\sigma_1(\overline{a\mbvdjord
b})(d)=\frac{1}{2}(a\mbvd
b+b\mbdv a)\mbvd d=\\
\frac{1}{2}(b\mbvd a\mbvd d+a\mbvd b\mbvd
d)=\frac{1}{2}(\sigma_1(\bar{a\vphantom b})\sigma_1(\bar
b)+\sigma_1(\bar b)\sigma_1(\bar{a\vphantom b}))(d).
\end{gather*}

(We write a composition of mappings as $fg(x)=g(f(x))$.) Analogously,
one may check that $\sigma_2$ is an associative specialization.

The relation (\ref{eq:SpecBiModCond}) follows from the definition
of the operation in our bimodule.

Moreover, $\sigma_1$ and $\sigma_2$ are commuting because
$$[\sigma_1(\bar{a\vphantom b}),\sigma_2(\bar b)](d)
=(\sigma_1(\bar{a\vphantom b})\sigma_2(\bar b) +\sigma_2(\bar
b)\sigma_1(\bar{a\vphantom b}))(d)=(a\mbvd d)\mbdv b-a\mbvd (d\mbdv
b)=0.$$

So, $J$ is a special $\bar J$-bimodule and by Theorem
\ref{thm:SpecSplitNullExtCrit} we obtain that $\widehat J$ is
special.
\end{proof}

In papers \cite{Kol:08,Kol:06} conformal algebras were investigated
and the following fact was proved.

\begin{prop}\label{prop:VarCur0DiVar}
If an algebra $A$ belongs to a variety $\Var$ then a dialgebra
$(\Cur A)^{(0)}$ belongs to a variety $\Di\Var$.
\end{prop}

We first prove an auxiliary statement.

\begin{lemma}\label{lemma:ifHomSJthenHomDiSJ}
If $\widehat J\in\Hom\SJ$ then $J\in\Hom\Di\SJ$.
\end{lemma}
\begin{proof}
Use conformal algebras. Let the algebra $\widehat J$ generated by a
set $X$ belong to the variety $\Hom\SJ$. Since $\SJ\,\langle
X\rangle$ is a free algebra of the variety $\Hom\SJ$, there exists
a surjective homomorphism $\phi\colon\SJ\,\langle X\rangle\to
\widehat J$. Then $\Cur\phi\colon\Cur\SJ\,\langle
X\rangle\to\Cur\widehat J$ is a morphism of conformal algebras and
particularly dialgebras. It is known \cite{GubKol:09} that
$J\hookrightarrow(\Cur\widehat J)^{(0)}$. So $(\Cur\phi)^{-1}[J]$
is a subdialgebra in $(\Cur\SJ\,\langle X\rangle)^{(0)}$. The
algebra $\SJ\,\langle X\rangle\in\SJ$ so by the definition of $\SJ$
there exists an associative algebra $A$ such that $\SJ\,\langle
X\rangle\hookrightarrow A^{(+)}$, hence $\Cur\SJ\,\langle
X\rangle\hookrightarrow\Cur A^{(+)}$ and $(\Cur\SJ\,\langle
X\rangle)^{(0)}\in\Di\SJ$. To complete the proof we need to
note that $J=\Cur\phi((\Cur\phi)^{-1}[J])$, where
$(\Cur\phi)^{-1}[J]\hookrightarrow(\Cur\SJ\,\langle
X\rangle)^{0}\in\Di\SJ$ and so $J\in\Hom\Di\SJ$.
\end{proof}

Now we can prove the following theorem.

\begin{thm}\label{thm:EqOfVarDialg}
$\Hom\Di\SJ=\Di\Hom\SJ$.
\end{thm}
\begin{proof}
To prove the inclusion "$\subseteq$" consider a free algebra
$\Di\SJ\,\langle X\rangle$ in the variety $\Hom\Di\SJ$. By Lemma
\ref{lemma:SpecUnitEmb} we have $\Di\SJ\langle
X\rangle\hookrightarrow J_e$, $J_e$ is special and $1\in \bar{J_e}$.
Then by Lemma \ref{lemma:SpecFact} $\bar J_e$ is special, hence
by Lemma \ref{lemma:SpecSplitNullExt} $\widehat{J_e}$ is a
special Jordan algebra and $J_e\in \Di\Hom\SJ$. Therefore,
$\Di\SJ\langle X\rangle\in\Di\Hom\SJ$. Since the free algebra of
the variety $\Hom\Di\SJ$ belongs to the variety $\Di\Hom\SJ$, the
variety $\Hom\Di\SJ$ is embedded into $\Di\Hom\SJ$.

We prove the inclusion "$\supseteq$". Let $J\in\Di\Hom\SJ$. By the
definition of a variety of dialgebras in the sense of Eilenberg it
means that $\widehat J\in\Hom\SJ$, hence by Lemma
\ref{lemma:ifHomSJthenHomDiSJ} we obtain $J\in\Hom\Di\SJ$.
\end{proof}

\subsection{s-identities in dialgebras}
Let $\Var$ be a variety of algebras, $X=\{x_1, x_2, \ldots \}$ be a
countable set. Consider a mapping $\phi_\Var\colon\Alg\,\langle X
\rangle \to \Var\,\langle X \rangle$ which maps $x_i\mapsto x_i$.
Let $T_0(\Var)$ be a set of multilinear polynomials from
$\ker\phi_\Var$, these are exactly all multilinear identities of
$\Var$. We suppose that the variety is defined by multilinear
identities that is $\Var=\{A\mid A\vDash T_0(\Var)\}$. There we use
the denotation $A\vDash f$ which means that the identity $f(x_1,
\ldots, x_n)=0$ holds on the algebra $A$.

Further, let $\Di\Alg 0\,\langle X\rangle$ be a free 0-dialgebra,
$\phi_{\Di\Var}\colon\Di\Alg 0\,\langle X \rangle \to
\Di\Var\,\langle X \rangle$, $T_0(\Di\Var)$ be a set of multilinear
dipolynomials from $\ker \phi_{\Di\Var}$, i.~e., all multilinear
identities from $\Di\Var$.

In paper \cite{Pozh:09} the following theorem was proved.

\begin{thm}[{Pozhidaev \cite[Theorem 3.2]{Pozh:09}}]\label{thm:DefDiVar}
Let $D\in\Di\Alg 0$. Then the following conditions are equivalent:
\begin{enumerate}
  \item $D\in\Di\Var$\textup{;}
  \item $\widehat D=\bar D\oplus D\in\Var$ \textup{(}the definition
  in the sense of Eilenberg\textup{);}
  \item $D\vDash \Psi^{x_i}_\Alg\,f$ for every $f\in T_0(\Var)$,
  $\deg f=n$, $i=1,\,\ldots,\,n$ \textup{(}the definition in the sense
  of \textup{\cite{Kol:08}).}
\end{enumerate}
\end{thm}

We prove the following
\begin{prop}\label{prop:PsiDiVarVar}
Let $f=f(x_1,\ldots,x_n)\in \Di\Alg 0\,\langle X\rangle$ be
multilinear, $f=\Psi^{x_j}_\Alg\,\bar f$ for some $j$. Then
$$f\in T_0(\Di\Var)\Leftrightarrow \bar f\in T_0(\Var).$$
\end{prop}
\begin{proof}
Since evidently $\Var\subseteq\Di\Var$, the statement
"$\Rightarrow$" is trivial.

To prove "$\Leftarrow$" consider an identity $\bar f\in T_0(\Var)$. By
Theorem \ref{thm:DefDiVar} for arbitrary $D\in\Di\Var$ we have
$D\vDash \Psi^{x_i}_\Alg\,\bar{f}$ for all $i=1,\,\ldots,\,n$, but
$\Psi^{x_j}_\Alg\,\bar{f}=f$ and so $f\in T_0(\Di\Var)$.
\end{proof}

\begin{prop}\label{prop:f1fnDiVarDiVar}
Let $f=f(x_1,\ldots,x_n)\in\Di\Alg 0\,\langle X\rangle$ be
multilinear, $f=f_1+\ldots+f_n$ where $f_i$ consists of all
dimonomials in $f$ with a central letter $x_i$. Then
$$f\in T_0(\Di\Var)\Leftrightarrow f_i\in T_0(\Di\Var)\text{ for all }
i=1,\,\ldots,\,n.$$
\end{prop}
\begin{proof}
"$\Leftarrow$" is evident.

We prove "$\Rightarrow$". Let $f\in T_0(\Di\Var)$, consider an
arbitrary algebra $A\in\Var$. Then by Proposition
\ref{prop:VarCur0DiVar} we obtain ${(\Cur A)}^{(0)}\in
\Di\Var$, hence ${(\Cur A)}^{(0)}\vDash f$, where
$f=f(x_1,\ldots,x_n)$. Fix $i\in\{1,\,\ldots,\,n\}$ and assign the following values to
variables: $x_i:=Ta_i$, $a_i\in A$, $x_j:=a_j$
for all $j\not=i$, $a_j\in A$. The properties of a conformal
product imply
$$0=f(a_1,\ldots,Ta_i,\ldots,a_n)=T\bar{f_i}(a_1,\ldots,a_n).$$

From the last equality we obtain $\bar{f_i}(a_1,\ldots,a_n)=0$,
so $A\vDash\bar{f_i}$ and $\bar{f_i}\in T_0(\Var)$. By the previous
proposition $f_i\in T_0(\Di\Var)$.
\end{proof}

We recall that $f$ is called a multilinear s-identity (in the case of
ordinary algebras) if
$$f\in T_0(\Hom\SJ)\setminus T_0(\Jord):=\SId.$$

A similar notion can be introduced for dialgebras \cite{Br:09}
$$f\in T_0(\Hom\Di\SJ)\setminus T_0(\Di\Jord):=\Di\SId.$$

\begin{thm}[about the correspondence of multilinear s-identities]\label{thm:CorrespSId}
\begin{enumerate}
  \item Let $g = g(x_1,\ldots,x_n)\in\SId$. Then $\Psi^{x_i}_\Alg\,g\in\Di\SId$ for all $i=1,\ldots,n$.
  \item Let $f=f(x_1,\ldots,x_n)\in\Di\SId$, $f=f_1+\ldots+f_n$
  \textup{(}by a central letter\textup{)}. Then there exists $j\in\{1,\,\ldots,\,n\}$
  such that $\bar{f_j}\in\SId$.
\end{enumerate}
\end{thm}
\begin{proof}
We prove the statement 1. Let $g\in\SId$, hence by the definition
$\SId$ we have $g\in T_0(\Hom\SJ)$ and $g\not\in T_0(\Jord)$. Proposition
\ref{prop:PsiDiVarVar} implies
$\Psi^{x_i}_\Alg\,g\in T_0(\Di\Hom\SJ)$, $\Psi^{x_i}_\Alg\,g\not\in
T_0(\Di\Jord)$. It follows from the equality of varieties
$\Hom\Di\SJ=\Di\Hom\SJ$ that $\Psi^{x_i}_\Alg\,g\in \Di\SId$.

For proving the statement 2 consider  $f\in\Di\SId$. By the
definition of $\Di\SId$ and Theorem \ref{thm:EqOfVarDialg} we
have $f\in T_0(\Hom\Di\SJ)=T_0(\Di\Hom\SJ)$ and $f\not\in
T_0(\Di\Jord)$. It follows from $f\in T_0(\Di\Hom\SJ)$ by Proposition
\ref{prop:f1fnDiVarDiVar} that $f_i\in T_0(\Di\Hom\SJ)$ for all $i$.
It follows from $f\not\in T_0(\Di\Jord)$
that $j$ exists such that $f_j\not\in T_0(\Di\Jord)$. Further, by
Proposition \ref{prop:PsiDiVarVar}, $\bar f_i\in T_0(\Hom\SJ)$
and $\bar f_j\not\in T_0(\Di\Jord)$, hence by the definition $\SId$
we obtain $\bar f_j\in\SId$.
\end{proof}

Now we can easily obtain the following corollary which was proved in
\cite{Br:09} by computer algebra methods.

\begin{cor}
There are no s-identities for Jordan dialgebras of degree
$\le 7$ and there exists a multilinear s-identity of a degree
8.
\end{cor}
\begin{proof}
Let $f$ be a s-identity for Jordan dialgebras, $\deg f=k\le 7$.
After a full linearization of $f$ we can suppose that $f$ is
multilinear that is $f\in\Di\SId$ and $f=f_1+\ldots+f_k$ by central letters.
It follows from Theorem \ref{thm:CorrespSId} about the
corresponding of multilinear s-identities that $\bar f_i\in\SId$ for some
$i$, $\deg \bar f_i\le k$, but Glennie proved \cite{Glennie:70} that
such an identity does not exist.

It is known \cite{Glennie:66} that there exists $f$ which is a
s-identity for Jordan algebras, $\deg f=8$. Again we can suppose
that $f$ is multilinear. Then Theorem \ref{thm:CorrespSId}
implies $\Psi^{x_i}_\Alg f$ is a required multilinear
s-identity for all $i=1,\,\ldots,\,8$.
\end{proof}

\subsection{Analogues for dialgebras of Shirshov's and Macdonald's Theorems}
Since we get the generalization of the Cohn's Theorem to the case of
dialgebras, a question appears about a generalization of the
Shirshov's Theorem for special Jordan algebras (\cite{Zhevl:78}, the simplification
of Shirshov's original proof is contained in \cite{JacPage:57}): Whenever every Jordan dialgebra
with two generators is special? The
answer to this questions is negative, it follows from Theorem
\ref{thm:Example2GeneratedExceptionalDialgebra}. However, the
following analogue of the Shirshov's Theorem holds for dialgebras.

\begin{thm}
Let $J$ be a one-generated dialgebra. Then $J$ is special.
\end{thm}
\begin{proof}
We have $J\in\Di\Jord$. Then by the definition a variety of
dialgebras in the sense of Eilenberg $\bar J\in\Jord$, $\widehat
J=\bar J\oplus J\in\Jord$. Let $x$ be the generator of $J$.
Then $\widehat J=\langle \bar x,x\rangle$, so $\widehat J$ is a
two-generated Jordan algebra. By the Shirshov's Theorem we obtain
that $\widehat J$ is special. We have $J\hookrightarrow(\Cur\widehat
J)^{(0)}$ and so $J$ is special too.
\end{proof}

Consider the particular case when two-generated dialgebra is free.

\begin{thm}
Let $J=\Di\Jord\,\langle x,y\rangle$ be the free Jordan dialgebra
generated by $x$, $y$. Then $J$ is special.
\end{thm}
\begin{proof}
We need to show that $J\in\Di\SJ$. First, prove
$J\in\Hom\Di\SJ$. Assume the converse, i.~e., $J\not\in\Hom\Di\SJ$.
By Lemma \ref{lemma:ifHomSJthenHomDiSJ} we obtain $\widehat
J=\bar J\oplus J\not\in\Hom\SJ$. Since $\widehat
J\in\Jord\setminus\Hom\SJ$, there exists a multilinear s-identity
$f(x_1,\ldots,x_n)$ of Jordan algebras such that $\SJ\vDash
f$ but $\widehat J\nvDash f$. Therefore, we may find
$u_1,\ldots,u_n\in\widehat J$ such that
$f(u_1,\ldots,u_n)\not =0$. Since the polynomial $f$ is
multilinear, we can suppose that either $u_i\in\bar J$ or $u_i\in J$ for
all $i$. The number of elements $u_i\in J$ does not exceed one otherwise,
$f(u_1,\ldots,u_n)=0$ because $J\cdot J=0$. Consider two
possible cases. The first case is when all $u_i\in\bar J$.
Then $\bar J\nvDash f$, which is impossible since $\bar J\in\SJ$
and $f$ is an s-identity. The second case is when
$u_1,\,\ldots,\,u_{n-1}\in\bar J$, $u_n\in J$. The algebra $\bar J$ is
generated by $\bar x$ and $\bar y$, so $u_i=u_i(\bar x,\bar y)$, $i=1,\ldots,n$. Then denote
$g(\bar x,\bar y,u_n):=f(u_1(\bar x,\bar y),\ldots,u_{n-1}(\bar
x,\bar y),u_n)\not=0$. Note that $g$ does not hold on $\widehat J$.
The polynomial $g(x,y,z)$ vanishes in $\SJ$, $\deg_z g=1$, hence by the
Macdonald's Theorem we obtain $g=0$ in $\Jord$. The
contradiction obtained proves that $J\in\Hom\Di\SJ$.

We prove that $J\in\Di\SJ$. We know that $J$ is a homomorphic
image of some special Jordan algebra $J_0$ under some mapping
$\phi\colon J_0\to J$. Let $x_0$ and $y_0$ are preimages of $x$
and $y$ with respect to $\phi$. Consider a
subdialgebra $U$ in $J_0$ generated by $x_0$ and $y_0$. Since the
dialgebra $J_0$ is special, subdialgebra $U$ is special too. The
dialgebra $J=\Di\Jord\,\langle x,y\rangle$ is free in the variety of
Jordan dialgebras, hence every mapping of $x$ and $y$ to
$U$ extends to a homomorphism. Map $x$ and $y$ to $x_0$ and $y_0$
respectively. Since $x_0$ and $y_0$ generate $U$, we obtain a
surjective homomorphism inverse to a homomorphism
$\phi|_U$. Therefore, $J\backsimeq U$ is a special Jordan dialgebra.
\end{proof}

\begin{cor}
If an identity $f(x,y)$ in two variables holds in all special
Jordan dialgebras then it holds in all Jordan dialgebras.
\end{cor}
\begin{proof}
Consider $f(x,y)$ as an element of the free Jordan dialgebra
$\Di\Jord\,\langle x,y\rangle$. By the previous theorem
$\Di\Jord\,\langle x,y\rangle$ is a special Jordan dialgebra,
therefore ${\Di\Jord\vDash f}$.
\end{proof}

In the paper \cite{Br:09} the s-identity of dialgebras was found
which depends on three variables and is linear in one of variables.
So the naive generalization of the Macdonald's Theorem to the
case of dialgebras is not true. But if an identity is linear in the
central letter then the following theorem is true which is an
analogue of the Macdonald's Theorem.

\begin{thm} Let $f=f(x,y,\dot z)$ be a dipolynomial which is linear
in $z$. If $\Di\SJ\vDash f$ then $\Di\Jord\vDash f$.
\end{thm}
\begin{proof}
Let $\Di\SJ\vDash f$, then $\Hom\Di\SJ\vDash f$. Consider a Jordan algebra $\bar
J\in\Hom\SJ$ as a dialgebra $J$ with equal
left and right products. Then $\bar J\in\Hom\SJ$ and $\widehat J=\bar
J\oplus J=\bar J\oplus\bar J\in\Hom\SJ$, so
$J\in\Di\Hom\SJ=\Hom\Di\SJ$. We obtain $J\vDash f$, hence $\bar
J\vDash \bar f$. Therefore, $\Hom\SJ\vDash\bar f=f(x,y,z)$, so by the
classical Macdonald's Theorem we have $\Jord\vDash\bar f$.

It remains to note that if $f= f(x,y,\dot z)$ is a
multilinear dipolynomial such that $\Jord\vDash\bar f=f(x,y,z)$
then $\Di\Jord\vDash f$: It follows immediately from the definition
\cite{Kol:08} of what is a variety of dialgebras. The
polynomial $f(x,y,z)$ can be nonlinear in $x$ and $y$. Suppose
$\deg_x f=n$, $\deg_y f=m$. Consider the full linearization
$$g(x_1,\dots, x_n, y_1,\dots, y_m, z)= L_x^n L_y^m f(x,y,z)$$ of
the identity $f(x,y,z)$ (notations from \cite[ch.~1]{Zhevl:78}).
Then $\Jord\vDash g(x_1,\dots, x_n, y_1,\dots, y_m, z)$ and so
$\Di\Jord \vDash g(x_1,\dots, x_n,y_1,\dots, y_m, \dot z)$.

If we now identify variables, then $$g(x,\dots, x, y,\dots, y,\dot
z)=n!m!f(x,y,\dot z).$$ In this section the characteristic of the
basic field is equal to zero, so we can divide by $n!m!$ and hence
$f(x,y,\dot z)$ is an identity on $\Di\Jord$.
\end{proof}

\begin{rem}
P.~M.~Cohn in \cite{Cohn:59} proposed an axiomatic characterization of
Jordan algebras $J_1(A)=A^{(+)}$ and $J_2(A,*)=H(A,*)$, where $A$ is an
associative algebra and $*$ is an involution on $A$, in terms of
$n$-ary operations. This is an interesting task to generalize these
results to Jordan dialgebras.
\end{rem}

\subsection*{Acknowledgements}
In the end of paper the author thanks P.~S.~Kolesnikov, A.~P.~Pozhidaev
and V.~Yu.~Gubarev for helpful discussions and valuable comments.
The author is grateful to the referee for valuable comments that allowed
to improve the manuscript. In particular, the statement about tetrads
in Theorem \ref{thm:CohnForDialgebra}, the criterion in
Proposition \ref{prop:CriterionOfQuotientSpecialityTwoGenerated}, and the final remark
were proposed to the author by the referee.

\end{document}